\documentclass[12pt,a4paper]{amsart}
\advance\hoffset-15truemm

\nonstopmode

\usepackage{paralist}
\usepackage{amssymb}
\usepackage{graphicx}

\newtheorem{thm}{Theorem}[section]

\newtheorem*{prop*}{Proposition}
\newtheorem*{thm*}{Theorem}

\newtheorem*{lem*}{Lemma}

\newtheorem*{cor*}{Corollary}

\theoremstyle{definition}

\newtheorem*{defn*}{Definition}
\theoremstyle{remark}
\newtheorem{rem}[thm]{Remark}
\newtheorem*{remk*}{Remark}

\newlength{\equwidth}
\def\A{\mathcal A}

\def\X{\mathfrak X}

\def\ga{\gamma}
\def\de{\delta}

\def\eps{\varepsilon}

\def\la{\lambda}
\def\rh{\rho}

\def\si{\sigma}
\def\ta{\tau}

\def\Ups{\Upsilon}
\def\ph{\varphi}

\def\Ga{\Gamma}

\def\Th{\Theta}
\def\La{\Lambda}

\def\Ph{\Phi}

\def\pa{\partial}
\def\t{\otimes}

\def\del{\pa}
\def\delstar{\pa^*}

\def\goesto{\rightarrow}

\def\cinf{\ensuremath{\mathrm{C}^\infty}}

\def\bdot{\mathrm{\bullet}}
\def\na{\mathrm{\nabla}}
\def\cov{\mathrm{\nabla}}

\def\dcov{\mathrm{d^{\nabla}}}

\def\dtcov{\mathrm{d^{\tilde\nabla}}}
\def\lapl{\mathrm{\square}}
\def\triang{\mathrm{\bigtriangleup}}
\def\Lapl{\mathrm{\triang}}

\def\ker{\mathrm{ker}\ }
\def\im{\mathrm{im}\ }

\def\gr{\mathrm{gr}}

\def\CC{\ensuremath{\mathcal{C}}}
\def\ZZ{\ensuremath{\mathcal{Z}}}
\def\BB{\ensuremath{\mathcal{B}}}
\def\HH{\ensuremath{\mathcal{H}}}

\def\ST{\ensuremath{\mathcal{S}}}

\def\co{\ensuremath{\mathfrak{co}}}

\def\so{\ensuremath{\mathfrak{so}}}

\def\l{\ensuremath{\mathfrak{l}}}
\def\a{\ensuremath{\mathfrak{a}}}
\def\g{\ensuremath{\mathfrak{g}}}

\def\gl{\ensuremath{\mathfrak{gl}}}

\def\ce{\ensuremath{\mathcal{E}}}
\newcommand{\bg}{\mbox{\boldmath{$ g$}}}

\def\Hol{\mathrm{Hol}}

\newcommand{\lpl}{
  \mbox{$
  \begin{picture}(12.7,8)(-.5,-1)
  \put(2,0.2){$+$}
  \put(6.2,2.8){\oval(8,8)[l]}
  \end{picture}$}}

\date{\today}
\title[Invariant Prolongation in Conformal Geometry]{Invariant Prolongation of BGG-Operators in Conformal Geometry}
\author{Matthias Hammerl}
\address{Institut f\"ur Mathematik, Universit\"at Wien, Nordbergstra\ss e~15, A--1090 Wien, Austria}
\email{matthias.hammerl@univie.ac.at}


\thanks{The author's work was supported by the IK I008-N funded by the University of Vienna}

\begin{document}

\begin{abstract}
BGG-operators form sequences of
invariant differential operators and 
the first of these is overdetermined. Interesting equations
in conformal geometry described by these operators are those for Einstein scales, conformal
Killing forms and conformal Killing tensors.
We present a deformation procedure of the tractor connection which
yields an invariant prolongation of the first operator. The explicit
calculation is presented in the case of conformal Killing forms.
\end{abstract}

\maketitle

\subsubsection*{\textbf{Acknowledgments}}
Special thanks go to Andreas \v{C}ap for his support and encouragement.
Discussions with Josef \v{S}ilhan have been very valuable and we are
indebted to his cogent remarks on an early version of this paper.
We thank Katja Sagerschnig for her comments on a draft. Discussions
with Vladimir Sou\v{c}ek have been very informative and we thank him for suggesting improvements for the notation.
The referee's careful criticism and suggestions for further remarks were
much appreciated.

\section{Introduction: Geometric prolongation of overdetermined operators. }\label{sec1}
A conformal structure on a manifold $M$\ is an equivalence class
$[g]$\ of pseudo-Riemannian metrics, where two metrics
$g$\ and $\hat g$\ are equivalent iff there is a function $f\in\cinf(M)$\
such that $\hat g=e^{2f}g$.
The simplest way to explain what a conformally invariant operator
is, is to give an example:
regard the operator
\begin{align}\label{Einsteinop}
  \Th^g:\cinf(M)&\goesto S^2_0 T^*M,\\
  \si&\mapsto \bigl(DD \si + \si P\bigr)_{0}.
\end{align}
Here $D$\ is the Levi-Civita connection of a metric $g$\
in the conformal class,
$P=P_{ab}$\ is the Schouten-tensor, which is a trace-modification of the Ricci tensor,
and the subscript $0$\ takes the trace-free part.
$S^2_0 T^*M$\ denotes
symmetric, trace-free bilinear forms on $TM$, which will
also be written as $\ce_{(ab)_0}$: Throughout the paper we
are using Penrose's abstract index notation \cite{penrose-rindler-87},
with $\ce^a=\Ga(TM)$\
denoting vector fields, $\ce_a=\Ga(T^*M)$\ denoting $1$-forms and
multiple indices being tensor products. Square and round brackets
around indices will indicate alternation resp. symmetrization.

Now $\Th^g$\ describes the equation governing Einstein scales:
for $\si\in\cinf(M)$\ one has $\Th^g \si=0$\ iff $\si^{-2}g$\
is Einstein.
The operator $\Th^g$\ is conformally covariant between
$\cinf(M)$\ and $S^2_0 T^*M$: if one switches to another
metric $\hat g=e^{2f}g$\ in the conformal class,
then
\begin{align*}
  {\Th}^{\hat g}\circ m(e^f)=m(e^f)\circ \Th^g,
\end{align*}
where $m(e^f)$\ is simply the multiplication operator
with $e^f$.
This yields a conformally invariant operator
between the weighted bundles $\HH_0=\ce[1]$\
and $\HH_1=S^2_0 T^*M\t\ce[1]$. Here
$\ce[w]$\ is the bundle of conformal $w$-densities,
which is a line bundle that is trivialized by every
metric $g$\ in the conformal class; with
$[\si]_g$\ the trivialization of a section $\si\in \ce[w]$,
one has $[\si]_{\hat g}=e^{wf}[\si]_g$.
Especially, $[g]$\ gives rise to a well defined \emph{conformal\ metric}\
$\bg=\bg_{ab}\in\ce_{(ab)}[2]$.

In general, a conformally invariant operator is obtained by a universal
formula in the Levi-Civita connection, the metric and the curvature,
possibly followed by contractions, that gives a well-defined operator
between natural bundles for the conformal structure.

The example of the operator for Einstein scales above has
another interesting property: it is overdetermined, and thus
one can wish to have a prolongation of the system: in classical
terms, this means that one wants to introduce more dependent
variables and derive differential consequences of the overdetermined
system, such that one can write down a closed system of equations;
i.e., a system of first order PDEs in which all (first order) derivatives
of the dependent variables are expressed in the dependent variables themselves.

\subsection{The standard tractor bundle of conformal geometry and the prolongation
of the equation governing Einstein scales}.

The prolongation of \eqref{Einsteinop}\ is well known, and is conformally
invariant. We are going to describe this and the necessary background on
conformal tractor bundles. Our notations are inspired by \cite{gover-silhan-2006}.
We note here that a reader who looks for an introduction to tractor calculus
in conformal geometry and an explanation of related notational issues could for instance
make use of the very careful and detailed exposition in the first part of \cite{josef-thesis}.

With respect to a metric $g$\ in the
conformal class the \emph{standard\ tractor\ bundle}\ $\ST$\ of a conformal
geometry is given by
\begin{align}\label{STg}
  [\ST]_g=\ce[1]\oplus \ce_a[1]\oplus\ce[-1]
\end{align}
and one writes elements $[s]_g=\si\oplus\ph_a\oplus\rh\in[\ST]_g$\ as
\begin{align}\label{column}
  [s]_g=
    \begin{pmatrix}
        \rh \\
        \ph_a \\
        \si
    \end{pmatrix}.
\end{align}
We remark here that via the conformal metric $\bg_{ab}\in\ce_{(ab)}[2]$\
and its inverse $\bg^{ab}\in\ce^{(ab)}[-2]$
one can move indices up and down, and thus we can also write
$[\ST]_g=\ce[1]\oplus \ce^a[-1]\oplus\ce[-1]$.

For $\hat g=\exp^{2f} g$\ one has the transformation
\begin{align}\label{transformS}
        [s]_{\hat g}=
        \begin{pmatrix}
            \hat\rh \\
            \hat\ph_a \\
            \hat\si
        \end{pmatrix}=
        \begin{pmatrix}
            \rh-\Ups_a\ph^a-\frac{1}{2}\si\Ups^b\Ups_b \\
            \ph_a+\si \Ups_a  \\
            \si
        \end{pmatrix}
\end{align}
where $\Ups=df$\ 
and $\ST$\ is defined by the equivalence class of $[\ST]_g$\ for
$g\in [g]$\ with respect to this transformation (\cite{thomass}).
We see that we have a well defined semi-direct composition series
\begin{align}\label{STsemi}
  \ST=\ce[1]\lpl\ce_a[1]\lpl\ce[-1]
\end{align}
 i.e., $\ST$\
is filtered $\ST=\ST^{-1}\supset\ST^{0}\supset\ST^{1}$
and with respect to a metric $g$\ in the conformal class
this filtration splits according to \eqref{STg}.

Additionally, $[\ST]_g$\ is endowed with the connection
\begin{align}\label{sttracon}
  \na_c s=\na_c
  \begin{pmatrix}
    \rh \\
    \ph_a \\
    \si
  \end{pmatrix}
  =
  \begin{pmatrix}
    D_c \rh-P_{c}^{\; b}\ph_b \\
    D_c\ph_a+\si P_{ca}+\rh \bg_{c a}\\
    D_c \si-\ph_c
  \end{pmatrix},
\end{align}
which is invariant with respect to the transformation \eqref{transformS}\
and thus gives a well defined connection on $\ST$, called the standard
tractor connection.

We furthermore see from \eqref{transformS}\ that one has a well-defined
projection $\Pi$\ to the `lowest slot' $\HH_0$\ of $\ST$.
This projection splits via the differential operator $L:\HH_0\goesto\ST$,
which is again defined via a metric $g$:
\begin{align*}
  \si\in\ce[1]\mapsto
  \begin{pmatrix}
    -\frac{1}{n}(\Lapl \si+{P_a}^{a}\si) \\
    \na\si \\
    \si
  \end{pmatrix}.
\end{align*}

$(\ST,\na,\Pi,L)$\ is a \emph{geometric\ prolongation}\ of $\Th:\HH_0\goesto \HH_1$:
The maps $\Pi$\ and $L$\ restrict to inverse isomorphisms of the
space of parallel sections of $\ST$\ with respect to $\na$\
and the space of Einstein scales in $\HH_0$.
This is well known. In the following section \ref{sec2}\ we will
give an explanation of this fact in terms of the BGG-machinery
and present a method to obtain invariant geometric prolongations
for other equations. In section \ref{sec3}\ we will give
an explicit prolongation of conformal Killing-forms \eqref{confkillequ}\ via this method.

\section{Conformal tractor bundles}\label{sec2}

We note here that the standard tractor bundle $\ST$\ and its tractor connection,
introduced via a description with respect to metrics in the conformal class above,
can alternatively be described as the associated bundle to the standard representation of
the Cartan-group $SO(p+1,q+1)$\ modelling conformal structures. More generally (see \cite{cap-gover-tractor},\cite{cap-gover-irr_tractor}): a tractor
bundle comes about as the associated space to a $SO(p+1,q+1)$-representation
and is canonically endowed with its tractor connection.

Apart from spin representations, all tractor bundles appearing in conformal
geometry appear as subbundles in tensorial powers of $\ST$, and we will assume
here for convenience of presentation that our tractor bundles are of this form.

Given $\ST$\ as in the previous section, i.e., written in terms
of a Levi-Civita connection $g$, one has a natural inner product
of signature $(p+1,q+1)$, which is given by
\begin{align}\label{inprod}
    \left(
      \begin{matrix}
      0 & 0 & 1 \\
      0 & g & 0 \\
      1 & 0 & 0       
      \end{matrix}
    \right).
\end{align}
Especially, one can identify $\La^2\ST$\ with $\so(\ST)$,
which is the \emph{adjoint}\ tractor bundle $\A M$\ for conformal structures.
Employing matrix notation, we will write its elements, or sections, as
\begin{align}\label{ammatrix}
  \begin{pmatrix}
    -c & -\eta_b & 0 \\
    \xi^a & C & \eta^b \\
    0 & -\xi_a & c
  \end{pmatrix},
\end{align}
where $(c,C)\in\co(p,q)$, $\xi^a\in\ce^a$\ and $\eta_a\in\ce_a$.

One has a natural surjection (projection to $\xi^a$) of $\A M$\ onto $TM$\
and an inclusion (inserting of $\eta_b$) of $T^*M$\ into $\A M$, 
while the inclusion via \eqref{ammatrix}\ of $TM$\ depends on the choice of $g$.
Having fixed a metric $g$\ in the conformal class, 
the algebraic action $\bdot$\ of $\A M$\ on a tractor bundle $T$\ restricts
to actions of $TM$\ and $T^*M$. Therefore, regarding $TM$\ and $T^*M$\ as (pointwise)
abelian Lie algebras, we can thus introduce Lie algebra differentials on
the the spaces $\CC_k:=\ce_{[c_1\cdots c_k]}\t T$:
we define $\del:\CC_k\goesto\CC_{k+1}$,
\begin{align}\label{formuladel}
  \del\ph(\xi_0,\cdots,\xi_k)&=\sum\limits_{j=0}^{k}(-1)^j\xi_j\bdot\ph(\xi_0,\cdots,
   \hat{\xi_j},\cdots,\xi_k)
\end{align}
and $\delstar:\CC_{k+1}\goesto\CC_{k}$,
\begin{align}\label{formuladelstar}
  \delstar Z_0\wedge\cdots\wedge Z_k\t V&=\sum\limits_{j=0}^k(-1)^{j+1}
  Z_0\wedge\cdots\wedge \widehat{Z_j}\cdots\wedge Z_k \t (Z_i\bdot V).
\end{align}

It is straightforward to check that $\del\circ\del=\delstar\circ\delstar=0$.
It is a consequence of a general result by Kostant (\cite{kostant-61}), that $\del$\ and $\delstar$\
are naturally adjoint with respect to
an (pointwise) inner product on the chain spaces $\CC_k$. This gives
a Hodge decomposition
\begin{align}
  \label{hodge}
  \CC_k=\im\del\oplus\ker\lapl\oplus\im\delstar
\end{align}
with $\lapl=\del\circ\delstar+\delstar\circ\del$.

Only $\delstar$, but not $\del$, is invariant with respect to a change
in Levi-Civita connection in the conformal class. Thus we use $\delstar$\
to define the spaces $\ZZ_k=\ker\delstar\cap \CC_k$,$\BB_k=\im\delstar\cap \CC_k$\
and $\HH_k=\ZZ_k/\BB_k$. Using the Hodge decomposition, one can identify
$\HH_k$\ with $\ker\lapl\subset \CC_k$\ after choice of a metric in the
conformal class.

As an $\so(\ST)$-valued form $K\in\ce_{[c_1c_2]}\t \A M$, the curvature
of the standard tractor connection is
\begin{align}\label{formulaK}
  K_{c_1c_2}=
  \begin{pmatrix}
    0 & -A_{ec_1c_2} & 0 \\
    0 & C_{c_1c2\;\; d}^{\;\ \;\; \ c} & A^e_{\;\ c_1c_2} \\
    0 & 0 & 0
  \end{pmatrix}.
\end{align}
Here $C$\ is the Weyl curvature and
$A=A_{ec_1c_2}=2D_{[c_1}P_{c_2]e}$\ is the Cotton-York tensor.
We recall that both $A_{ec_1c_2}$\ and $C_{c_1c2\quad d}^{\;\ \;\ c}$\ are trace-free. Furthermore
the skew-symmetrization over any $3$ indices of
$C_{abcd}$\ vanishes,
as does the skew-symmetrization of $A_{abc}$.
The Weyl curvature doesn't satisfy the differential Bianchi
identity, however one has
\begin{align*}
  D_{[a} C_{bc]de}=\bg_{d[a}A_{|e|bc]}-\bg_{e[a}A_{|d|bc]}.
\end{align*}

\section{Invariant geometric prolongation via the BGG-machinery}\label{sec3}

The BGG-machinery will associate to the tractor covariant derivative $\na$\ on $T$\ 
differential operators $\Th_l:\HH_l\goesto \HH_{l+1}$.
The core step in the construction of the BGG-operators is to find, in a natural way,
a splitting of $\Pi_l:\ZZ_l\goesto \HH_l$\ adapted to $\na$:
For every $\si\in \HH_l$\ it can be shown that there is
a unique lift $s\in \ZZ_l$\ such that $\dcov s\in \ZZ_{l+1}$.
This defines the \emph{BGG-splitting\ operators} $L_l:\HH_l\goesto
\ZZ_{l}$. By construction they give rise to the \emph{BGG-operators}
\begin{align*}
    \Th_l&:\HH_l\goesto \HH_{l+1},\\
    \Th_k&=\Pi_{k+1}\circ \dcov\circ L_k.
\end{align*}
Thus one obtains the BGG-sequence
\begin{align*}
  0 \goesto \HH_0\overset{\Th_0}{\goesto} \HH_1 {\goesto}\cdots
\overset{\Th_{n-1}}{\goesto} \HH_{n}\goesto 0.
\end{align*}

We are interested in the first operator $\Th_0$, which gives
an overdetermined system of equations. 
In \cite{prolong}\ a prolongation
method for operators $\Th_0+\eta$\ for $\eta$\ a lower order, possibly nonlinear,
differential operator, was developed, which did not however take into
account the invariance respectively naturality of the operator $\Th_0$.

The construction of BGG-operators sketched above also works for more general
connections $\tilde\na=\na+\Psi$, if $\Psi\in\ce_c\t\gl(T)^1$. Here
$\gl(T)^1$\ denotes those endomorphisms of $T$\ which are homogeneous
of degree $\geq 1$\ with respect to the filtration of $T$\ inherited from
$\ST$. More simply put:\ $\gl(T)^1$\ consists of upper triangular matrices if
we use vector-notation as in \eqref{column}\ and \eqref{Vup}.

\subsection{Deformation of the tractor connection}

We would like to understand the solution space of $\Th_0:\HH_0\goesto \HH_1$:
Let $\si\in \HH_0$. By definition
  $\Th_0\si=\Pi_1(\na(L_0\si)),$
and thus
\begin{align}\label{nullspace}
  \Th_0\si=0\ \mathrm{iff}\ \na(L_0\si)\in \BB_1=\ker \Pi_1,
\end{align}
which shows that in general $(T,\na,\Pi_0,L_0)$\ is
not a prolongation, since the kernel of $\Th_0$\ is not
mapped into the space of parallel sections by $L_0$: while
parallel sections of $(T,\na)$\ always project into
the kernel of $\Th_0$,  a solution $\si\in\HH_0$\
of $\Th_0\si=0$\ will, by definition of $\Th_0$, only
have the property that $\na(L_0\si)\in\im\delstar$.
Our strategy is to deform $\na$\ to $\tilde\na=\na+\Psi$\
by a map $\Psi\in\ce_c\t\gl(T)^1$\ in suitable way,
such that we obtain a `better' on $T$\ which
gives a geometric prolongation of $\Th_0$. I.e.,
we want to find $\tilde\na$\ such that $\Th_0\si=0$\
implies $\tilde\na(L_0\si)=0$\ and conversely.

We make the following observation: consider $\Psi\in\ce_c\t\gl(T)^1$\
which has the property that
\begin{align}
  \label{normprop1}
  \Psi s\in \im\delstar=\BB_1\ \forall s\in T.
\end{align}
Then we can construct the BGG-splitting operators $\tilde L_0:H_0\goesto T$,
$\tilde L_1:H_1\goesto \CC_1$\ and the first BGG-operator
$\tilde \Th_0:\HH_0\goesto\HH_1$\ as above. But Since $(\tilde\na-\na)=\Psi$\
has values in $\BB_1$, we see that 
$\delstar\circ \tilde \nabla\circ L_0=\delstar\circ\nabla\circ L_0=0$,
which shows that we have $\tilde L_0=L_0$; and since $\Pi_1(\BB_1)=0$,
we have $\tilde \Th_0=\Th_0$. Thus: maps $\Psi\in\ce_c\t\gl(T)^1$\
which send $T$\ into $\BB_1$\ may be used to deform $\na$\ to another
connection without changing the first BGG-operator.
Thus the space of such $\Psi$\ gives us a freedom for suitable deformations
of $\na$.

Assume that we have managed to find such a $\Psi$\ for which
the curvature $R$\ of $\tilde\na=\na+\Psi$\ has the property
that, for every $s\in\Ga(T)$,
\begin{align}\label{curveprop}
   \delstar(Rs)=0.
\end{align}
Then we claim that $\tilde \na s= \tilde L_1 \Th_0 \Pi_0(s)$.
This means that for every $s\in \Ga(T)$, one already has
$\delstar(\dtcov(\tilde\na s))$. But this expression
equals $\delstar(Rs)$, and thus we have the claimed commutativity.

But this is already enough: because now, if $\Th_0\si=0$,
we have that $\tilde\na(L_0\si)=\tilde L_1(0)=0$. And on the other
hand, for a parallel section $s$\ of $T$, one evidently
has by construction of $L_0$\ that $L_0(\Pi_0(s))=s$.
Thus, $\Pi_0:T\goesto \HH_0$\ and $L_0:\HH_0\goesto T$\ restrict
to inverse isomorphisms between parallel sections of $T$\ with respect
to $\tilde\na$\ and the kernel of $\Th_0$.

Therefore the whole problem lies in finding a deformation map
$\Psi\in\ce_c\t \gl(T)^1$\ which maps $T$\ into $\im\delstar$\
and which gives a $\tilde\na=\na+\Psi$\ whose curvature $R$\
maps $T$\ into $\ker\delstar\subset\ce_{[c_1c_2]}\t T$.
Existence and uniqueness of such a map can be shown using
analogs of inductive normalization procedures well known in the
realm of parabolic geometries and this prolongation method actually works in a
more general situation. See also remark \ref{remaweyl}. 

One constructs the new connection $\tilde\na$\ in terms
of a given metric $g$\ in the conformal class.
By uniqueness of $\Psi$, this construction is however
independent of the choice of $g$, i.e., the connection $\tilde\na$\
is a well defined, conformally invariant object.
In the following we are going to show how this construction
of a deformation map $\Psi$\ works explicitly
for a special and interesting case,
that of conformal Killing forms.

\section{Invariant prolongation of conformal Killing forms.}\label{sec4}

Conformal Killing forms were first prolonged by U. Semmelmann \cite{semmelmann}, however the discussion there did not take into account conformal invariance
of the equation. In \cite{gover-silhan-2006}\ an invariant prolongation
was calculated by ad hoc methods (see also \cite{josef-thesis}). The following is a completely conceptual
derivation of an invariant geometric prolongation by describing conformal
Killing forms as the kernel of a first BGG-operator and prolonging
this operator via a deformation of the tractor connection as in section
\ref{sec3}.

We are going to proceed as follows: In 
\ref{subsec1}\ we describe the exterior powers
of the standard tractor bundles, give explicit formulas for the Lie algebraic
differentials on the first chain spaces and determine their
$CO(p,q)$-decompositions.
In \ref{subsec2}\ we describe explicitly how the operator governing
conformal Killing $k$-forms comes about as first BGG-operator for
the $k+1$-st exterior power of the standard tractor bundle. 
In \ref{subsec3}\ we obtain a geometric prolongation by
constructing a deformation $\Psi\in\ce_c\t \gl(T)^1$\ with
the properties called for in section \ref{sec3}.
In \ref{subsec4}\ we show the conformal invariance of $\Psi$.

\subsection{The tractor bundle.}
\label{subsec1}
In the following $k$\ will be $\geq 1$. The tractor bundle 
$T=\La^{k+1}\ST$\ decomposes (via a metric $g$\ in the conformal class) into
\begin{align}\label{Vup}
    \begin{pmatrix}
        \ce_{[a_1\cdots a_k]}[k-1] \\
        \ce_{[a_1\cdots a_{k+1}]}[k+1]\; |\;\ \ce_{[a_1\cdots a_{k-1}]}[k-1] \\
        \ce_{[a_1\cdots a_k]}[k+1]
    \end{pmatrix}
\end{align}
and similarly as for \eqref{transformS}, we have transformations
\begin{align}\label{transformT}
        \begin{pmatrix}
            \hat\rh_{a_1\cdots a_k} \\
            \hat\ph_{a_0\cdots a_k}\; |\;\ \hat\mu_{a_2\cdots a_k} \\
            \hat\si_{a_1\cdots a_k}
        \end{pmatrix}=
        \begin{pmatrix}
            \rh_{a_1\cdots a_k}-\Ups^b\ph_{b a_1\cdots a_k}
                              -k\Ups_{[a_1}\mu_{a_2\cdots a_k]}
                            -\frac{1}{2}\Ups^b\Ups_b \si_{a_1\cdots a_k} \\
            \ph_{a_0\cdots a_k}+(k+1)\Ups_{[a_0}\si_{a_1\cdots a_k]} \; |\; 
            \mu_{a_2\cdots a_k}-\Ups_b\si_{b a_2\cdots a_k}
            \\
            \si_{a_1\cdots a_k}
        \end{pmatrix}.
\end{align}
From \eqref{STsemi}, or directly from \eqref{transformT}, we see that
\begin{align}\label{Tsemi}
  T=\ce_{[a_1\cdots a_k]}[k+1]\lpl(\ce_{[a_1\cdots a_{k+1}]}[k+1]
\oplus\ce_{[a_1\cdots a_{k-1}]}[k-1])
\lpl\ce_{[a_1\cdots a_k]}[k-1],
\end{align}
which splits into $T_{-1}\oplus T_0\oplus T_1$\ after choice of $g$\ in the conformal class.

The standard tractor connection \eqref{sttracon}\ gives rise
to the invariantly defined tractor connection $\na$\
on $T$:
\begin{align}\label{tracon}
    &\na_c
    \begin{pmatrix}
            \rh_{a_1\cdots a_k} \\
            \ph_{a_0\cdots a_k}\; |\;\ \mu_{a_2\cdots a_k} \\
            \si_{a_1\cdots a_k}
    \end{pmatrix}  
    =
\left(
\begin{matrix}
      {D}_c\rh_{a_1\cdots a_k} -P_{c}^{\; p}\ph_{pa_1\cdots a_k}-kP_{c[a_1}\mu_{a_2\cdots a_k]} \\
\left(
\begin{matrix}
      {D}_c\ph_{a_0\cdots a_k} +(k+1)\bg_{c[a_0}\rh_{a_1\cdots a_k]} \\  
      +(k+1)P_{c[a_0}\si_{a_1\cdots a_k]}
\end{matrix}
\right)
      \; | \;\
\left(
\begin{matrix}
 {D}_c\mu_{a_2\cdots a_k} 
      \\
      -P_{c}^{\; p}\si_{pa_2\cdots a_k}
      +\rh_{ca_2\cdots a_k}
\end{matrix}
\right)
      \\
      {D}_c\si_{a_1\cdots a_k}
      - \ph_{ca_1\cdots a_k}+k\de_{c[a_1}\mu_{a_2\cdots a_k]}.
    \end{matrix}
    \right).
\end{align}

\subsubsection{Description of the first homology groups}\label{V-homology}
$\delstar:\CC_1\goesto \CC_0=T$\ is given (see \eqref{formuladelstar})
by $Z\t s\mapsto -Z\bdot s$\ for $s\in \Ga(T)$, $Z\in\ce_a$.
Thus $B_0=\im\delstar:\CC_1\goesto T$\ is
simply $\ce_a\bdot T$, which is all of $T^0$.
Thus $\HH_0=T/T^0$.
By the Hodge decomposition \eqref{hodge}\ we
can embed $\HH_0$\ as $T_{-1}=\ker\lapl=\ker\del\subset T$.

Also,
$\HH_i$\ will be embedded into $\CC_i$\ as $\ker\lapl=\ker(\del\delstar+\delstar\del)$\
for $i=1,2$.
The calculation of the $CO(p,q)$-decomposition of the spaces $\HH_i$
is purely algorithmic using Kostant's
version of the Bott-Borel-Weyl theorem \cite{kostant-61};
the details of which are not important for us here.
We just state the results for $\HH_1$\ and $\HH_2$, which are all homologies
we are going to need:
We will write
\begin{align*}
  \CC_i=
  \begin{pmatrix}
    \ce_{[c_1\cdots c_k]}\t T_1 \\
    \ce_{[c_1\cdots c_k]}\t T_0 \\
    \ce_{[c_1\cdots c_k]}\t T_{-1}
  \end{pmatrix},
\end{align*}
and speak of the top, middle and bottom slots.

$\ce_c\t T_{-1}$\ contains the highest weight part
$\ce_{\{c[a_1\ldots a_k]\}_0}[k+1]$, and this
is all of $\HH_1$.
Explicitly,
$\ce_{\{c[a_1\ldots a_k]\}_0}[k+1]$\
sits in $\ce_{c[a_1\ldots a_k]}[k+1]$
as those $\si=\si_{ca_1\cdots a_k}$\
which have both zero trace and vanishing alternation:
\begin{align*}
  0=g^{pq}\si_{pq a_2\cdots a_k},\quad
  0=\si_{[c a_2\cdots a_k]}.
\end{align*}

If $k\geq 2$\ then the analogous statement holds also for the
second chain space: in this case $\HH_2$\ is exactly the highest
weight part of $\ce_{[c_1c_2]}\t T_{-1}=\ce_{[c_1c_2][a_1\ldots a_k]}[k+1]$.
i.e.,
$\HH_2=\ce_{\{[c_1c_2][a_1\ldots a_k]\}_0}[k+1]\subset\ce_{[c_1c_2]}\t T_{-1}$.

Especially, for $i=0,1$\ we have that $H_i$\ lies in
the lowest grading part of $\CC_i$\ and if $k\geq 2$\ this
also holds for $i=2$:
\begin{align*}
  \begin{pmatrix}
    T_1 \\
    T_0   \\
    \HH_0=T_{-1}
  \end{pmatrix}
  \overset{\del}{\goesto}
  \begin{pmatrix}
    \ce_c\t T_1 \\
    \ce_c\t T_0   \\
    \HH_1\oplus \im\del_{|T_0}
  \end{pmatrix}
  \overset{\del}{\goesto}
  \begin{pmatrix}
    \ce_{[c_1c_2]}\t T_1 \\
    \ce_{[c_1c_2]}\t T_0   \\
    \HH_2\oplus \im\del_{|\ce_c\t T_0}
  \end{pmatrix}
\end{align*}

Now we describe what $\del$,$\delstar$\ and $\lapl$\ do
on the first few chain spaces $\CC_0=T,\CC_1=\ce_c\t T$\
and $\CC_2=\ce_{[c_1c_2]}\t T$:

\subsubsection{Explicit formulas for $\del,\delstar$\ and $\lapl$\ on
  the first chain spaces.}\label{B-explicit}

\begin{align}\label{delformulas}
        \del_c
        \begin{pmatrix}
            \rh_{a_1\cdots a_k} \\
            \ph_{a_0\cdots a_k}\ |\   \mu_{a_2\cdots a_k} \\
            \si_{a_1\cdots a_k}
        \end{pmatrix}
      &=
      \begin{pmatrix}
        0 \\
       (k+1)\de_{c[a_0}\rh_{a_1\cdots a_k]}  \ |\
       \rh_{ca_2\cdots a_k} \\
       -\ph_{ca_1\cdots a_k} + k\bg_{c[a_1}\mu_{a_2\cdots a_k]}
      \end{pmatrix} \\ \notag
        \del_{c_1}
        \begin{pmatrix}
            \rh_{c_2 a_1\cdots a_k} \\
            \ph_{c_2 a_0\cdots a_k}\ |\   \mu_{c_2 a_2\cdots a_k} \\
            \si_{c_2 a_1\cdots a_k}
        \end{pmatrix}
      &=
      \begin{pmatrix}
        0 \\
        2(k+1)\bg_{[a_0|[c_1}\rh_{c_2]|a_1\cdots a_k]}
        \ |\ 
        -2\rh_{[c_1c_2] a_2\cdots a_k} \\
        2\ph_{[c_1c_2] a_2\cdots a_k} 
        +2k\bg_{[a_1|[c_1}\mu_{c_2]| a_2\cdots a_k]}
      \end{pmatrix}\\ \notag
        \del^*
        \begin{pmatrix}
            \rh_{c a_1\cdots a_k} \\
            \ph_{c a_0\cdots a_k}\ |\   \mu_{c a_2\cdots a_k} \\
            \si_{c\ a_1\cdots a_k}
        \end{pmatrix}
      &=
      \begin{pmatrix}
        \bg^{pq}\ph_{p q a_1\cdots a_k} + k\mu_{[a_1\cdots a_k]} \\
        -(k+1)\si_{[a_0\cdots a_k]}
        \ |\ 
        \bg^{pq}\si_{p qa_2\cdots a_k} \\
        0
      \end{pmatrix}\\ \notag
        \del^*
        \begin{pmatrix}
            \rh_{c_1c_2 a_1\cdots a_k} \\
            \ph_{c_1c_2 a_0\cdots a_k}\ 
            |\   \mu_{c_1c_2 a_2\cdots a_k} \\
            \si_{c_1c_2 a_1\cdots a_k}
        \end{pmatrix}
      &=
      \begin{pmatrix}
        -2\bg^{pq}\ph_{cpqa_1\cdots a_k} -2k\mu_{c[a_1\cdots a_k]} \\
        2(k+1)\si_{c[a_0\cdots a_k]}
        \ |\
        -2\bg^{pq}\si_{cp q a_2\cdot a_k} \\
        0
      \end{pmatrix}.
\end{align}
The image of $\delstar$\ in $T=C_0$\ is
simply $T^0=T_0\oplus T_1$, and the Kostant Laplacian
thus acts by positive real scalars on $T_1$\ and the two
components of $T_0$. It
vanishes on $T_{-1}$\ by \eqref{hodge}. Explicitly, $\lapl$
is given on $T$\ by
\begin{align}\label{kost1}
  \begin{pmatrix}
    n \\
    (k+1)
    \ |\
    (n-k+1) \\ 
    0
  \end{pmatrix}.
\end{align}

The image of $\delstar$\ in $C_1$\ contains all
of $\ce_c\t T_1$\ (since we have \eqref{hodge}). Now $\ce_c\t T_1$\ decomposes
into three parts: the alternating maps, $\ce_{k+1}[k+1]$,
the purely trace maps, $\ce_{k-1}[k-1]$, and finally those
maps which have both trivial trace and trivial alternating part,
$\ce_{\{c[a_1\ldots a_k]\}}[k+1]$. We will denote the three
irreducible components of $\ce_c\t T_1$\ by $(\ce_c\t T_1)_{alt}$,
$(\ce_c\t T_1)_{\{\}_0}$\ and $(\ce_c\t T_1)_{tr}$.
We will write this decomposition of $\ce_c\t T_1\cap\im\delstar=\ce_c\t T_1$\
\begin{align}\label{decomp2}
  \begin{pmatrix}
    alt \\
    \{\}_0 \\
    tr
  \end{pmatrix},
\end{align}
and in this picture the Kostant Laplacian $\lapl$\ acts by
\begin{align}\label{kost2up}
  \begin{pmatrix}
    2(n+k-1) \\
    2(n-2) \\
    2(2n-k-1)
  \end{pmatrix}.
\end{align}

Now to the middle slot:
We have
\begin{align*}
  \ce_c\t T_0=\ce_{c [a_0\cdots a_k]}[k+1]\oplus\ce_{c [a_2\cdots a_k]}[k-1]
\end{align*}
and both parts split into alternating, $\{\}_0$-\ and trace components.
Both $\{\}_0$-components,
the left alternating and the right trace component lie in the image of $\delstar$.
The only other component of $\im\delstar\cap \ce_c\t T_0$\ is
$\ce_{[a_1\cdots a_k]}[k-1]$, which embeds into $\ce_c\t T_0$\ via
\begin{align*}
  \ta_{a_1\cdots a_k}\mapsto
  \begin{pmatrix}
    0 \\
    -k(k+1)\bg_{c[a_0}\ta_{a_1\cdots a_k}
    \ |\ 
    (n-k)\ta_{c a_2\cdots a_k}\\
    0
  \end{pmatrix}.
\end{align*}
We will write the decomposition of $\ce_c\t T_0\cap\im\delstar \subset \ce_c\t T_0$
\begin{align}\label{decomp1}
  \begin{pmatrix}
    alt &|\ &* \\
    \{\}_0 &|\ &\{\}_0 \\
    tr &|\ &tr
  \end{pmatrix},
\end{align}
and the Kostant Laplacian is seen to act by the scalars
\begin{align}\label{kost2mid}
  \begin{pmatrix}
   4(k+1) &|\ &* \\
    2k &|\ &2(n-k) \\
    2n &|\ &2(n-k-1)
  \end{pmatrix}.
\end{align}

\subsection{The first BGG-operator  $\Th_0:\HH_0\goesto \HH_1$\
and conformal Killing forms.}
\label{subsec2}

Using $\eqref{tracon},\eqref{delformulas}$\ and $\eqref{kost1}$,
we compute that the first BGG-splitting operator
$L_0:\HH_0\goesto T$\ is, up to first homogeneity, given by
\begin{align}\label{L0form}
  \si\mapsto
  \begin{pmatrix}
    *
    \\
    {D}_{[a_0}\si_{a_1\cdots a_k]} \; | \;
    -\frac{1}{n-k+1}\bg^{pq}{D}_p\si_{qa_2\cdots a_k}
    \\
    \si_{a_1\cdots a_k}
  \end{pmatrix}.
\end{align}

In \ref{V-homology}\ we saw that (using $\bg\in\ce^{(ab)}[2]$),
\begin{align*}
    \HH_0&=\ce_{[a_1\cdots a_k]}[k+1],\\
    \HH_1&=\ce_{{\{c[a_1\cdots a_k]\}}_0}[k+1],\\
    \HH_2&=\ce_{{\{[c_1c_2][a_1\cdots a_k]\}}_0}[k+1].
\end{align*}

Thus we immediately see using \eqref{tracon}\ that
\begin{align*}
    \Th_0&:\HH_0\goesto\HH_1,\\
    \Th_0&=\Pi_1\circ\cov\circ L
\end{align*}
is given by
${\Th_0}_c \si_{a_1\cdots a_k}= {D}_{\{c}\si_{a_1\cdots a_k\}_0}$.
But this is exactly the conformal Killing operator.
Thus our prolongation procedure will yield an
isomorphism between the space conformal Killing forms
\begin{align}\label{confkillequ}
  \si_{a_1\cdots a_k}\in \ce_{[a_1\cdots a_k]}[k+1],\
  {D}_{\{c}\si_{a_1\cdots a_k\}_0}=0
\end{align}
and the space of parallel sections of a natural connection on $T$.

For $k=1$\ we get exactly the operator describing conformal Killing fields,
i.e., infinitesimal automorphisms of the conformal structure; see also remark \ref{remak1}.
This case has been treated in detail in \cite{gover-lapl_einstein}.  The
main result of this text, an explicit geometric prolongation, will also
work for $k=1$. We only need $k \geq 2$\ for obtaining an algebraic obstruction
tensor which is described in subsection \ref{subsubcurv}. 

\subsection{The deformation of the tractor connection}
\label{subsec3}

We are now going to construct a $\Psi\in\ce_c\t \gl(T)^1$\
with the properties called for in section \ref{sec3}.

The calculations will be made more readable by providing beforehand
the mappings which will appear:
We will make use of the vector bundle maps
\begin{align*}
  L_i:\ce_{[a_1\cdots a_k]}[k+1]\goesto \ce_{c [a_1\cdots a_{k+1}]}[k+1]\quad, i=1,2,
\end{align*}
and
\begin{align*}
  R_i:\ce_{[a_1\cdots a_k]}[-k+1]\goesto \ce_{c [a_1\cdots a_{k-1}]}[k-1]\quad, i=1,2,
\end{align*}
of homogeneity $1$:
\begin{align}\label{maps1}
  L_1(\si)&=C_{[a_0a_1\ c}^{\quad\;\ p}\si_{|p|a_2\cdots a_k]} &
  L_2(\si)&=\bg_{c[a_0}C_{a_1a_2}^{\quad\;\ pq}\si_{|pq|a_3\cdots a_k]} \\ \notag
  R_1(\si)&=C_{c[a_2}^{\quad\; pq}\si_{|pq|a_3\cdots a_k]} &
  R_2(\si)&=C_{[a_2a_3}^{\quad\; pq}\si_{c|pq|a_4\ldots a_k]}.
\end{align}
In homogeneity $2$\ we will need the maps 
\begin{align*}
  F_i,G_i:\ce_{[a_1\cdots a_k]}[k+1]\goesto \ce_{c[a_1\cdots a_k]}[k-1],
\end{align*}
the maps 
\begin{align*}
  E_i:\ce_{[a_1\cdots a_{k+1}]}[k+1]\goesto \ce_{c[a_1\cdots a_k]}[k-1]
\end{align*}
and the maps 
\begin{align*}
  T_i:\ce_{[a_1\cdots a_{k-1}]}[k-1]\goesto \ce_{c[a_1\cdots a_k]}[k-1]:
\end{align*}
\begin{align}\label{maps2}
    E_1(\ph)&=C_{c[a_1}^{\quad\; pq}\ph_{|pq|a_2\cdots a_k]} &
    E_2(\ph)&=C_{[a_1a_2}^{\quad\;\ pq}\ph_{c|pq|a_3\cdots a_k]} \\ \notag
    T_1(\mu)&=C_{c\;\; [a_1a_2}^{\;\ p}\mu_{|p|a_3\cdots a_k]} &
    T_2(\mu)&=\bg_{c[a_1}C_{a_2a_3}^{\quad\;\ pq}\mu_{|pq|a_3\cdots a_k]}\\ \notag
    F_1(\si)&=A_{[a_1c}^{\quad\; p}\si_{|p|a_2\cdots a_k]} &
    F_2(\si)&=A^p_{\; c[a_1}\si_{|p|a_2\cdots a_k]}\\ \notag
    F_3(\si)&=A^p_{\; [a_1a_2}\si_{c|p|a_3\cdots \a_k]} &
    F_4(\si)&=\bg_{c[a_1}A_{a_2}^{\;\ pq}\si_{|pq|a_3\cdots\a_k]}\\ \notag
    G_1(\si)&=(D_cC_{[a_1a_2}^{\quad\; pq})\si_{|pq|a_3\cdots a_k]} &
    G_2(\si)&=(D^pC_{c\;[a_1a_2}^{\; q})\si_{|pq|a_3\cdots a_k]}\\ \notag
    G_3(\si)&=(D_{[a_1}C_{|c| a_2}^{\quad\; pq})\si_{|pq|a_3\cdots a_k]}.
\end{align}

With respect to the $CO(p,q)$-decompositions \eqref{decomp1}\ and
\eqref{decomp2}\ a more natural basis for the linear space formed
by these maps into of $\gr(\CC_1)_1$\ and $\gr(\CC_1)_2$ is formed by

\begin{align}\label{maps-basis}
    L_{tr} &=-\frac{k-1}{n-k}L_2 &
    L &=L_1-L_{tr} \\ \notag
    R_{alt} &=\frac{2}{k}R_1+\frac{k-2}{k}R_2 &
    R &=\frac{k-2}{k}(R_1-R_2) \\ \notag
    E_{alt} &=\frac{2}{k+1}E_1+\frac{k-1}{k+1}E_2 &
    E &=\frac{k-1}{k+1}(E_1-E_2) \\ \notag
    T_{tr} &=-\frac{k-2}{n-k+1}T_2 &
    T &=T_1-T_{tr} \\ \notag
    F_{tr} &=\frac{k}{n-k+1}F_4 &
    F_{alt} &=\frac{2}{k+1}F_2-\frac{k-1}{k+1}F_3 \\ \notag
    F_i &=F_1-\frac{1}{k+1}F_2+\frac{k-1}{2(k+1)}F_3 -\frac{k-1}{k}F_{tr} &
    F_{ii} &=\frac{k-1}{k+1}(F_2+F_3) -\frac{k-1}{2k}F_{tr} \\ \notag
    G_{i} &=G_1+2 F_{alt}-\frac{2}{k}(n-k-1)F_{tr} &
    G_{ii} &=G_2-\frac{2(k-2)}{k}F_{tr} \\ \notag
    G_{iii} &=G_3-2 F_{alt}-\frac{n-3}{k}F_{tr}.
\end{align}

$L_{tr},T_{tr}$\ and $F_{tr}$\ are purely trace,
$R_{alt}$,$E_{alt}$\ and $F_{alt}$\ are alternating and all
other maps have both vanihing alternation and trace.

The maps of \eqref{maps1}\ and \eqref{maps2}\ can be expressed as
\begin{align}\label{maps-decomp}
    L_1&=L+L_{tr} &
    L_2&=-\frac{n-k}{k-1}L_{tr} &
    R_1&=R+R_{alt} 
    \\ \notag
    R_2&=-\frac{2}{k-2}R+R_{alt} &
    E_1&=E+E_{alt} &
    E_2&=-\frac{2}{k-1}E+E_{alt}
    \\ \notag
    T_1&=T+T_{tr} &
    T_2&=-\frac{n-k+1}{k-2}T_{tr}
\end{align}
and
\begin{align*}
    F_1&=F_{i}+\frac{1}{2}F_{alt}+\frac{k-1}{k}F_{tr} &
    F_2&=F_{ii}+F_{alt}+\frac{k-1}{2k}F_{tr} \\ \notag
    F_3&=\frac{2}{k-1}F_{ii}-F_{alt}+\frac{1}{k}F_{tr} &
    F_4&=\frac{n-k+1}{k}F_{tr} \\ \notag
    G_1&=G_{i}-2F_{alt}+\frac{2}{k}(n-k-1)F_{tr} &
    G_2&=G_{ii}+\frac{2(k-2)}{k}F_{tr} \\ \notag
    G_3&=G_{iii}+2 F_{alt}+\frac{n-3}{k}F_{tr}.
\end{align*}

For
$  s=
        \begin{pmatrix}
            \rh_{a_1\cdots a_k} \\
            \ph_{a_0\cdots a_k}\; |\;\ \mu_{a_2\cdots a_k} \\
            \si_{a_1\cdots a_k}
        \end{pmatrix}
$
we have
\begin{align*}
  (K\bdot s)=
  \begin{pmatrix}    
    C_{c_1c_2 [a_1}^{\quad\quad p}\rh_{|p|a_2\cdots a_k]}
    -k A_{[a_1|c_1c_2|}\mu_{a_2\cdots a_k}
    -A^p_{\; c_1c_2}\ph_{|p|a_1\cdots a_k}
    \\
    C_{c_1c_2 a_0}^{\quad\quad p}\ph_{pa_1\cdots a_k]} 
    +A_{[a_0| c_1c_2|}\si_{a_1\cdots a_k]}
    \; |\; 
    C_{c_1c_2[a_2}^{\quad\quad p}\mu_{pa_3\cdots a_k]}
    -A^p_{\; c_1c_2}\si_{|p|a_2\cdots a_k}
    \\
    C_{c_1c_2 [a_1}^{\quad\quad p}\si_{pa_2\cdots a_k]}
  \end{pmatrix}.
\end{align*}
We calculate
\begin{align}\label{delstarkappa}
    \delstar(K\bdot s)=
    \begin{pmatrix}
      2k F_1 + 2k F_2 -k E_1 +k(k-1)T_1 \\
      -k(k+1)L_1\ |\ -(k-1)R_1 \\
      0
    \end{pmatrix}
\end{align}
and thus have that the lowest homogeneous component of $\delstar(K\bdot s)$,
which is of homogeneity $1$, is given by $(-k(k+1) L_1 |\  -(k-1)R_1)$.
Now we use \eqref{maps-decomp},\eqref{kost1}\ and $\eqref{maps-basis}$
to apply $-\lapl^{-1}$\ to this expression, which yields
\begin{align}\label{psi1def}
    \Psi_1:=
    \begin{pmatrix}
      0 \\
      \la_1 L_1+\la_2 L_2  |\ \rh_1 R_1+\rh_2 R_2 \\
      0
    \end{pmatrix}
\end{align}
where
\begin{align*}
  \la_1&=\frac{1+k}{2} & 
  \la_2&=\frac{(k-1)(k+1)}{2n} \\
  \rh_1&=\frac{(k-1)(n-2)}{2(n-k)n} &
  \rh_2&=\frac{2-3k+k^2}{2(k-n)n}.
\end{align*}

Now the curvature of the deformed connection $\na+\Psi_1$\ is
\begin{align*}
  R=K\bdot+\dcov\Psi_1+[\Psi_1,\Psi_1],
\end{align*}
but $[\Psi_1,\Psi_1]$\ 
obviously vanishes.
Let us calculate $R$: The only term which demands our
attention is $\dcov\Psi_1$. Take any
$  s=
  \begin{pmatrix}
    \rh_{a_1\cdots a_k} \\
    \ph_{a_0\cdots a_k} | \mu_{a_2\cdots a_k} \\
    \si_{a_1\cdots a_k}
  \end{pmatrix}
  \in\Ga(T)$.

Then for $\xi_1,\xi_2\in\X(M)$, we have, since $\Psi_1$\ is
a $1$-form on $M$\ with values in $\gl(T)$,
\begin{align}\label{calc1}
  &(\dcov\Psi_1)s(\xi_1,\xi_2)=\\&= \notag
\cov_{\xi_1}(\Psi_1(\xi_2)s)-\Psi_1(\xi_2)(\cov_{\xi_1}s)
-\cov_{\xi_2}(\Psi_1(\xi_1)s)+\Psi_1(\xi_1)(\cov_{\xi_2}s)
-\Psi_1([\xi_1,\xi_2])s.
\end{align}
We may expand \eqref{calc1}\ and write
\begin{align}\label{calc2}
  &(\dcov\Psi_1)s=
  &\left(\begin{matrix}
   *\\
   \left(\begin{matrix}
 {D}_{\xi_1}\bigl({\Psi_1}(\xi_2)\si\bigr)
 -{\Psi_1}({\xi_2})\bigl({D}_{\xi_1}\si\bigr)
 -{D}_{\xi_2}\bigl({\Psi_1}(\xi_1)\si\bigr)
 +{\Psi_1}({\xi_1})\bigl({D}_{\xi_2}\si\bigr)
\\
 -\Psi_1([\xi_1,\xi_2])\si
\\
  -{\Psi_1}(\xi_2)\del_{\xi_1}\ph
  +{\Psi_1}(\xi_1)\del_{\xi_2}\ph
  -{\Psi_1}(\xi_2)\del_{\xi_1}\mu
  +{\Psi_1}(\xi_1)\del_{\xi_2}\mu        
   \end{matrix}\right)
    \\
    \del_{\xi_1}\Psi_1(\xi_2)\si
    -\del_{\xi_2}\Psi_1(\xi_1)\si
  \end{matrix}\right),
\end{align}
where we don't care about the top component since
it vanishes after an application of $\delstar$.
The lowest component is simply
$\del(\Psi_1\si)=-\del\lapl^{-1}\delstar(K\bdot\si).$
Thus $\delstar(Rs)$\ lies in the top slot (i.e., in homogeneity $1$).
So our first deformation had the effect of moving the
expression $\delstar\circ R$\ one slot higher.
This can be repeated:
it is a straightforward calculation using the expression
in the middle component of \eqref{calc2} and, in that order,
\eqref{delstarkappa},
\eqref{maps-decomp},\eqref{kost2up}\ and $\eqref{maps-basis}$
to find, with $\phi:=-\lapl^{-1}\circ\delstar\circ R$,
\begin{align}\label{psi2def}
    \Psi=\Psi_2=\Psi_1-\phi=
    \begin{pmatrix}
      \begin{pmatrix}
        \eps_1 E_1+\eps_2 E_2+\tau_1 T_1+\tau_2 T_2 \\
        +\phi_1 F_1+\phi_2 F_2+\phi_3 F_3+\phi_4 F_4 \\
        +\ga_1 G_1+\ga_2 G_2+\ga_3 G_3
      \end{pmatrix}
       \\
      \la_1 L_1+\la_2 L_2  |\ \rh_1 R_1+\rh_2 R_2 \\
      0      
    \end{pmatrix}
\end{align}
where
\begin{align*}
  \eps_1&=\frac{k-1}{2(n-k)} &
  \eps_2&=\frac{(k-1)k}{2(k-n)n} \\
  \tau_1&=\frac{(k-1)(n(n-k+1)-2k)}{2(k-n)n} &
  \tau_2&=-\frac{(k-2)(k-1)}{2n} \\
  \phi_1&=-\frac{n+k-3}{n-2} &
  \phi_2&=\frac{1-k}{n} \\
  \phi_3&=\frac{(k-1)(n+k)}{2(k-n)n} &
  \phi_4&=\frac{(k-1)(2+k-2n)}{2(k-n)(n-2)} \\
  \ga_1&=-\frac{k-1}{2(n-2)n} &
  \ga_2&=\frac{k-1}{2(n-2)} \\
  \ga_3&=\frac{(k-1)k}{2(k-n)n}.
\end{align*}

Now the curvature $R'$\ of $\na+\Psi=\na+\Psi_1+\phi$\
is given by
\begin{align*}
  R+\dcov\phi+[\Psi_1,\phi]+[\phi,\Psi_1]+[\phi,\phi].
\end{align*}
One sees that for every $s\in\Ga(T)$,
$([\Psi_1,\phi]+[\phi,\Psi_1]+[\phi,\phi])s$\ has only
values in the top component and
we may therefore forget about this term when calculating
$\delstar (R's)$. As in the calculation \eqref{calc2},
we see that $(\dcov\phi) s$\ has only values in the middle and top slots
and the middle slot is given by $2\del_{[c_1}\phi_{c_2]}s$.
Therefore, by construction of $\phi$, we see that
$\delstar(Rs)$\ vanishes, and thus, via
the considerations of section \ref{sec3},
we have solved the prolongation problem for conformal Killing forms.

We have already remarked there that this solution must already
be conformally invariant by virtue of uniqueness, which is not difficult
to see, but to see what is going on
we are going to check independence of the choice of metric by hand in \ref{subsec4}.

\begin{rem}\label{remak1}
  For $k=1$, we have $T=\La^2 \ST=\so(\ST)=\A M$. Thus $\na+\Psi$\
  prolongs the first BGG-operator for the adjoint tractor connection in this case.
  But in \cite{cap-infinitaut}\ it was shown that for a parabolic
  geometry which is either $1$-graded or torsion free and which has the property
  that the first homology of the adjoint tractor bundle $\HH_1$\ is concentrated
  in lowest homogeneity, the corresponding first BGG-operator describes infinitesimal
  automorphisms of the structure - in our case conformal Killing fields - and is
  prolonged by the connection
  $\tilde\na s=\na s+i_s K$\ 
  which maps $T$\ into $\im\delstar$\ and satisfies \eqref{curveprop}.
  Thus uniqueness implies that $\Psi s=i_s K$, which can also be read off \eqref{psi2def}\
  directly.
\end{rem}

\begin{rem}
  The invariant connection prolonging the conformal Killing equation \eqref{confkillequ}\
  which was constructed in \cite{gover-silhan-2006}\  differs from our result $\Psi$
  as defined in \eqref{psi2def}\ since it can be checked to have nontrivial intersection
  with $\im\del$. Recall that our solution obeys the normalization conditions \eqref{normprop1}\
  and \eqref{curveprop}, but the first of these conditions implies that it has
  trivial intersection with $\im\del$\ since one has the Hodge decomposition \eqref{hodge}.

  If one wants to translate the solution \eqref{psi2def}\ into the
  notation used in \cite{gover-silhan-2006}\ one has to use
  the automorphism
  \begin{align*}
    \left(
    \begin{matrix}     
      \rh_{a_1\cdots a_k} \\
      \ph_{a_0\cdots a_k} | \mu_{a_2\cdots a_k} \\
      \si_{a_1\cdots a_k}
    \end{matrix})
    \right)
    \mapsto 
    \left(
    \begin{matrix}     
      (k+1)\rh_{a_1\cdots a_k} \\
      \ph_{a_0\cdots a_k} | -k(k+1)\mu_{a_2\cdots a_k} \\
      (k+1)\si_{a_1\cdots a_k}
    \end{matrix})
    \right)
  \end{align*}
  which transforms an element of
  \begin{align*}
    \ce_{[a_1\cdots a_k]}[-k+1]\lpl(\ce_{[a_1\cdots a_{k+1}]}[-k-1]\oplus\ce_{[a_1\cdots a_{k-1}]}[-k+1])\lpl\ce_{[a_1\cdots a_k]}[-k-1]
  \end{align*}
in our notation to the equivalent
  element in the notation of \cite{gover-silhan-2006}.
  Then $\Psi$\ as defined in \eqref{psi2def}\ has the following
  form with respect to the notations of Gover-\v{S}ilhan:
  \begin{align*}
    \Psi_c(
    \left(
    \begin{matrix}     
      \rh_{a_1\cdots a_k} \\
      \ph_{a_0\cdots a_k} | \mu_{a_2\cdot sa_k} \\
      \si_{a_1\cdots a_k}
    \end{matrix})
    \right)
    =
    \left(
    \begin{matrix}
      \left(\begin{matrix}
        (k+1)(\eps_1E_1+\eps_2E_2)\ph-\frac{1}{k}(\ta_1T_1+\ta_2T_2)\mu \\  
        +(\phi_1 F_1+\phi_2 F_2+\phi_3 F_3+\phi_4 F_4 \\
        +\ga_1 G_1+\ga_2 G_2+\ga_3 G_3)\si
      \end{matrix}\right)    \\
      \frac{1}{k+1}(\la_1L_1+\la_2L_2)\si | -k(\rh_1R_1 + \rh_2R_2)\si \\
      0 
    \end{matrix}
    \right)
    .
  \end{align*}
\end{rem}

\subsubsection{Algebraic obstruction tensors obtained via the curvature of the deformed connection}\label{subsubcurv}
Since $\Psi\in\ce_c\t\gl(T)^1$\ we know that the curvature $R\in\ce_{[c_1c_2]}(\gl(T))$\ of $\tilde\na=\na+\Psi$\
agrees with $K\bdot$\ in homogeneity $0$. But if $\si_{a_1\cdots a_k}\in\ce^{k}[k+1]$\ is
a conformal Killing $k$-form, then $L_0\si$\ is given by \eqref{L0form}; and thus
$0=\dtcov\tilde\cov s=Rs$\ agrees with $K\bdot L_0\si$\ in $\ce_{[c_1c_2]}\t T_{-1}$. But by
\eqref{formulaK}\ this is simply (minus) $C_{c_1 c_2\;\ [a_1}^{\quad\;\; p}\si_{|p|a_2\cdots a_k]}$.
For $k\geq 2$\ we have $\HH_2=\ce_{\{[c_1c_2][a_1\cdots a_k]\}_0}[k+1]$\ and projecting the previous expression to this space gives the conformally
invariant algebraic map
\begin{align*}
  \Ph:\si\mapsto C_{\{c_1 c_2\;\ [a_1}^{\quad\;\;\ p}\si_{|p|a_2\cdots a_k]\}_0}
\end{align*}
which vanishes on conformal Killing $k$-forms.
This obstruction has also been constructed T. Kashiwada in \cite{kashiwada}, U. Semmelmann in \cite{semmelmann}\ and recently by R. Gover and J. \v{S}ilhan in \cite{gover-silhan-2006}.
Our derivation is completely conceptual: the map is simply the composition
of the first two BGG-operators for the deformed connection $\tilde\na$: $\Ph=\tilde \Th_1\circ \Th_0$.
This evidently explains both conformal invariance of $\Ph$\ and why it vanishes
on conformal Killing forms. That $\Ph$\ is algebraic has the cohomological reason
that $\HH_2$\ is concentrated in lowest homogeneity.

\begin{rem}
Apart from the (trivial) cases of Einstein scales and twistor spinors
where one doesn't need any deformation and automatically has \eqref{curveprop},
the case of conformal Killing forms is the simplest situation in
which to explicitly compute the prolongation. The next interesting
case to treat will be conformal Killing tensors, for which, as far as we know,
there has not yet been given any prolongation, and which can
be treated similarly as the form case.
There the situation becomes more complicated however, since
the modelling representations $S^k\g$\ are $2k+1$-graded. This case could have
interesting relations to symmetries of the Laplacian (\cite{eastwood-laplacian}).
\end{rem}

\begin{rem}
The holonomy $\Hol(\tilde\na)$\ of the thus obtained prolongation connection
$\tilde\na$\ describes
the solution space of the operator $\Th_0$. In the case
of the standard tractor bundle and the spinor tractor bundle\ 
one has $\na=\tilde\na$\ (see section \ref{sec3})
and thus the solution space is governed by the conformal holonomy
of the structure, i.e., existence of Einstein scales and twistor spinors correspond to reductions of the conformal holonomy. 
In general, the existence of
non-trivial solutions of $\Th_0$\ doesn't imply
a holonomy reduction: for instance, full conformal holonomy doesn't
obstruct the existence of conformal Killing fields or conformal Killing
forms.

Because of \eqref{nullspace}, parallel sections of a tractor bundle
give special solutions to $\Th_0$.
In the case of conformal Killing Forms, those  coming from parallel sections
were called \emph{normal}\ conformal Killing forms by F. Leitner in \cite{leitner-normal}. 
This notion of \emph{normal}\ solutions of first BGG-operators makes
sense for every tractor bundle and they correspond to reductions
in conformal holonomy.
\end{rem}

\begin{rem}
  Using the tractor approach above for describing Einstein scales
as parallel sections, R. Gover and P. Nurowski \cite{gover-nurowski-2006}\ used
the curvature $R$\ of the standard tractor connection and its derivatives to obtain (under
a genericity condition on the Weyl curvature) a conformally invariant
system of tensors which provides a sharp obstruction
against the existence of Einstein scales.
For a general tractor bundle and $R$\ the curvature of
the prolongation connection, one can similarly build natural
systems of obstruction tensors, but it is not known whether these will be sharp.
\end{rem}

\subsection{Conformal invariance of $\Psi$}
\label{subsec4}
For this calculation we need some transformation formulas.
We will denote by $\hat D$\ the Levi-Civita connection
of the metric rescaled by $e^{2f}$. More generally, we will
denote by a hatted symbol the corresponding quantity calculated with respect
to the metric $\hat g$. With $\Ups=df$\ we have
\begin{align*}
  {{\hat D}}_u C_{abcd}&={D}_u C_{abcd}-2\Ups_u C_{abcd}
  -2 \Ups_{[a}C_{|u|b]cd} - 2\Ups_{[c}C_{|u|d]ab}\\
  &\quad+2(n-3)g_{u[a}A_{b]cd} + 2 (n-3)g_{u[c} A_{d]ab}\\
  \hat A_{abc}&=A_{abc}+\Ups^d C_{dabc}\\
\end{align*}
In the calculation the following transformation-maps
\begin{align*}
  H_i:\ce^{[a_1\cdots a_k]}[-k+1]\goesto\ce^{[a_1\cdots a_k]}[-k-1]
\end{align*}
will appear:
\begin{align*}
   H_1(\si)&=\Ups_c C_{[a_1a_2}^{\quad\;\; pq}\si_{|pq|a_3\cdots a_k]} &
   H_6(\si)&=\bg_{c[a_1}\Ups^u C_{a_2 a_3}^{\quad\;\; pq}\si_{|upq|a_4\cdots a_k]} \\
   H_2(\si)&=\Ups_c C_{[a_1a_2}^{\quad\;\; pq}\si_{|pq|\a_3\cdots a_k]} &
   H_7(\si)&=\Ups_d C_{[a_1 a_2}^{\quad\;\;\ dp} \si_{c|p|a_3\cdots a_k]} \\
   H_3(\si)&=\Ups^p C_{[a_1 a_2 c}^{\quad\quad q}\si_{|pq|a_3\cdots a_k]} &
   H_8(\si)&=\Ups^d C_{d[a_1c}^{\quad\;\;\ p}\si_{|p|a_2\cdots a_k]} \\
   H_4(\si)&=\bg_{c[a_1}\Ups^{d}C_{|d| a_2}^{\quad\;\; pq}\si_{|pq|a_3\cdots a_k]} &
   H_9(\si)&=\Ups_d C_{ca_1}^{\quad dp}\si_{|p|a_2\cdots a_k]} \\
   H_5(\si)&=\Ups_{[a_1}C_{a_2a_3}^{\quad\;\ pq}\si_{|cpq|a_4\cdots a_k]}.
\end{align*}

The maps \eqref{maps1}\ of homogeneity $1$\ are invariant with respect to the choice
of $g\in [g]$\ since the Weyl curvature is conformally invariant.
It is straightforward to calculate that the individual maps \eqref{maps2} transform like
\begin{align*}
  \hat E_1&=E_1+2H_9-(k-1)H_2 &
  \hat E_2&=E_2+H_1-2H_7-(k-2)H_5 \\
  \hat G_1&=G_1-2H_1-2H_2-2H_3+2H_4+2H_7 &
  \hat G_2&=G_2-H_1-H_2-H_3+H_7+2H_8 \\
  \hat G_3&=G_3 -H_1-H_2-H_3+H_4+2H_9,
\end{align*}
and
\begin{align*}
  \hat F_1&=F_1+H_8 &
  \hat F_2&=F_2+H_9 &
  \hat F_3&=F_3+H_7 \\
  \hat F_4&=F_4+H_4 &
  \hat T_1&=T_1-H_3 &
  \hat T_2&=T_2-H_6.
\end{align*}

Thus, if we switch to another metric $\hat g$\ respectively the corresponding
linear connection  ${\hat D}$\
and then calculate $\hat\Psi$\ using \eqref{psi2def},
the result differs from $\Psi$\ only in the top slot of $C_1$\
and it does so by
\begin{align}\label{trans1}
    &(\eps_2-2\ga_1-\ga_2-\ga_3) H_1 
    -\ta_2 H_6
    -((k-1)\eps_1-2\ga_1-\ga_2-\ga_3) H_2 \\ \notag
    &+(-2\eps_2+\phi_3+2\ga_1+\ga_2)H_7 
    -(\ta_1-2\ga_1-\ga_2-\ga_3)H_3 
    +(\phi_1+2\ga_2)H_8 \\ \notag
    &+(\phi_4+2\ga_1+\ga_3)H_4 
    +(2\eps_1+\phi_2+2\ga_3)H_9
    -(k-2)\eps_2 H_5.
\end{align}
On the other hand, if we calculate $\Psi$\ with respect to $g$\ and then transform the
expression via
$\hat\rh=\rh-\Ups_d\ph^{da_1\cdots a_k}-k\Ups^{[a_1}\mu^{a_2\cdots a_k]},$
the difference to $\Psi$\ also lies in homogeneity two and is
\begin{align}\label{trans2}
    -\la_2\frac{1}{k+1}H_1 -\rh_1 H_2+\frac{k-1}{k+1}\la_1H_3
    +\frac{2}{k+1}\la_2 H_4
    -k \rh_2 H_5 + \frac{k-1}{k+1}\la_2H_6 -\frac{2}{k+1}\la_1H_8.
\end{align}
Now it is straightforward to check that the expressions \eqref{trans1}\
and \eqref{trans2}\ coincide.
Thus $\Psi$\ is seen not to depend on the choice of the metric in the
conformal class used
to construct it. As we already remarked this is in fact a
consequence of the uniqueness property of $\Psi$\
stated at the end of section \ref{sec3}.

\begin{rem}\label{remaweyl}
The prolongation method of above works more generally: the construction
of the BGG-sequence works for arbitrary tractor bundles
over regular $k$-graded parabolic geometries (\cite{BGG-2001},\cite{BGG-Calderbank-Diemer}), and again
the first operator in this sequence is overdetermined,
and we ask for a natural prolongation. The analog of
a choice of metric in the conformal case is the choice
of a Weyl structure of the parabolic geometry (\cite{cap-gover-tractor}). The
homogeneity conditions become a bit more subtle, but
the basic principle of finding a natural deformation
of the tractor connection yielding a prolongation is the same
as presented in section \ref{sec3}.
This is the subject of a forthcoming joint paper with J. \v{S}ilhan, P. Somberg and V. Sou\v{c}ek.
\end{rem}


\begin{thebibliography}{10}

\bibitem{thomass}
T.~N. Bailey, M.~G. Eastwood, and A.~Rod Gover.
\newblock Thomas's structure bundle for conformal, projective and related
  structures.
\newblock {\em Rocky Mountain J. Math.}, 24(4):1191--1217, 1994.

\bibitem{prolong}
Thomas Branson, Andreas {\v{C}}ap, Michael Eastwood, and A.~Rod Gover.
\newblock {Prolongations of geometric overdetermined systems.}
\newblock {\em Int. J. Math.}, 17(6):641--664, 2006.

\bibitem{BGG-Calderbank-Diemer}
David~M.J. Calderbank and Tammo Diemer.
\newblock {Differential invariants and curved Bernstein-Gelfand-Gelfand
  sequences.}
\newblock {\em J. Reine Angew. Math.}, 537:67--103, 2001.

\bibitem{cap-infinitaut}
Andreas {\v{C}}ap.
\newblock Infinitesimal {A}utomorphisms and {D}eformations of {P}arabolic
  {G}eometries.
\newblock {\em J. Europ. Math. Soc.}
\newblock To appear.

\bibitem{cap-conformal-bgg}
Andreas {\v{C}ap}.
\newblock Overdetermined systems, conformal geometry, and the {BGG} complex.
\newblock In M.G. Eastwood and W.~Millor, editors, {\em "Symmetries and
  Overdetermined Systems of Partial Differential Equations", The IMA Volumes in
  Mathematics and its Applications}, volume 144, pages 1--25. Springer, 2008.

\bibitem{cap-gover-irr_tractor}
Andreas {\v{C}}ap and A.~Rod Gover.
\newblock Tractor bundles for irreducible parabolic geometries.
\newblock In {\em Global analysis and harmonic analysis (Marseille-Luminy,
  1999)}, volume~4 of {\em S\'emin. Congr.}, pages 129--154. Soc. Math. France,
  Paris, 2000.

\bibitem{cap-gover-tractor}
Andreas {\v{C}}ap and A.~Rod Gover.
\newblock Tractor calculi for parabolic geometries.
\newblock {\em Trans. Amer. Math. Soc.}, 354(4):1511--1548 (electronic), 2002.

\bibitem{BGG-2001}
Andreas {\v{C}}ap, Jan Slov\'{a}k, and Vladim\'{\i}r Sou\v{c}ek.
\newblock {Bernstein-Gelfand-Gelfand sequences.}
\newblock {\em Ann. of Math.}, 154(1):97--113, 2001.

\bibitem{eastwood-laplacian}
Michael Eastwood.
\newblock Higher symmetries of the {L}aplacian.
\newblock {\em Ann. of Math. (2)}, 161(3):1645--1665, 2005.

\bibitem{gover-lapl_einstein}
A.~Rod Gover.
\newblock Laplacian operators and {$Q$}-curvature on conformally {E}instein
  manifolds.
\newblock {\em Math. Ann.}, 336(2):311--334, 2006.

\bibitem{gover-nurowski-2006}
A.~Rod Gover and Pawe{\l} Nurowski.
\newblock {Obstructions to conformally Einstein metrics in $n$ dimensions.}
\newblock {\em J. Geom. Phys.}, 56(3):450--484, 2006.

\bibitem{gover-silhan-2006}
A.~Rod Gover and Josef \v{S}ilhan.
\newblock The conformal {K}illing equation on forms -- prolongations and
  applications.
\newblock {\em Diff. Geom. Appl.}
\newblock To appear.

\bibitem{kashiwada}
T.~Kashiwada.
\newblock {On conformal Killing tensor.}
\newblock {\em Natur. Sci. Rep. Ochanomizu Univ.}, 19:67--74, 1968.

\bibitem{kostant-61}
Bertram Kostant.
\newblock Lie algebra cohomology and the generalized {B}orel-{W}eil theorem.
\newblock {\em Ann. of Math. (2)}, 74:329--387, 1961.

\bibitem{leitner-normal}
Felipe Leitner.
\newblock Conformal {K}illing forms with normalisation condition.
\newblock {\em Rend. Circ. Mat. Palermo (2) Suppl.}, (75):279--292, 2005.

\bibitem{penrose-rindler-87}
Roger Penrose and Wolfgang Rindler.
\newblock {\em Spinors and space-time. {V}ol.\ 1}.
\newblock Cambridge Monographs on Mathematical Physics. Cambridge University
  Press, Cambridge, 1987.
\newblock Two-spinor calculus and relativistic fields.

\bibitem{semmelmann}
Uwe Semmelmann.
\newblock {Conformal Killing forms on Riemannian manifolds.}
\newblock {\em Math. Z.}, 245(3):503--527, 2003.

\bibitem{josef-thesis}
Josef \v{S}ilhan.
\newblock {\em Invariant operators in conformal geometry}.
\newblock PhD thesis, University of Auckland, 2006.
\newblock \\ \verb+http://www.math.muni.cz/~silhan/ps/Thesis_UoA.ps+.

\end{thebibliography}
\end{document}